\documentclass[12pt]{article}
\usepackage[T1]{fontenc}
\usepackage{amsthm}
\usepackage{amsmath,amsfonts,txfonts,amssymb,mathrsfs}
\usepackage{lmodern}
\usepackage[all]{xy}
\title{}
\author{Saurav Bhaumik}

\newtheorem{thm}{Theorem}[section]
\newtheorem{lem}[thm]{Lemma}
\newtheorem{proposition}[thm]{Proposition}
\newtheorem{defin}{Definition}
\newtheorem{corol}[thm]{Corollary}

\newenvironment{lemma}{\begin{lem}\;}{\end{lem}}

\newenvironment{theorem}{\begin{thm}\;}{\end{thm}}

\theoremstyle{remark}
\newtheorem{remark}[thm]{\bf Remark}

\mathchardef\mhyphen="2D

\newcommand{\gm}{\mathbb G_m}

\newcommand{\ftwo}{\mathbb F_2}

\newcommand{\inv}{{}^{-1}}
\newcommand{\isomor}{\,\raisebox{4pt}{$\sim$}{\kern -.89em\to}\,}

\newcommand{\spec}{{\rm Spec\;}}

\newcommand{\topet}{{}_{\rm\acute{e}t}}
\newcommand{\topliset}{{}_{\rm lis\mhyphen\acute{e}t}}
\newcommand{\hliset}{\!\!\!\topliset}
\newcommand{\het}{\!\!\topet}
\newcommand{\CC}{\mathbb C}

\newcommand{\tr}{T_R}
\newcommand{\ta}{T_A}

\newcommand{\hsing}{\!{}_{\rm sing}}
\title{The Chow ring for the classifying space of $GO(2n)$}
\author{Saurav Bhaumik}
\date{}
\begin{document}
\maketitle
\begin{abstract}
Let $GO(2n)$ be the general orthogonal group scheme (the group of orthogonal similitudes). In the topological category, Y. Holla and N. Nitsure determined the singular cohomology ring $H^*_{\rm sing}(BGO(2n,\mathbb C),\mathbb F_2)$ of the classifying space $BGO(2n,\mathbb C)$ of the corresponding complex Lie group $GO(2n,\mathbb C)$ in terms of explicit generators and relations. The author of the present note showed that over any algebraically closed field of characteristic not equal to $2$, the smooth-\'etale cohomology ring $H_{\rm sm\mhyphen\'et}^*(BGO(2n),\mathbb F_2)$ of the classifying algebraic stack $BGO(2n)$ has the same description in terms of generators and relations as the singular cohomology ring $H^*_{\rm sing}(BGO(2n,\mathbb C),\mathbb F_2)$. Totaro defined for any reductive group $G$ over a field, the Chow ring $A^*_G$, which is canonically identified with the ring of characteristic classes in the sense of intersection theory, for principal $G$-bundles, locally trivial in \'etale topology. In this paper, we calculate the Chow group $A^*_{GO(2n)}$ over any field of characteristic different from $2$ in terms of generators and relations.
\end{abstract}
\section{Introduction}
The Chow ring of the classifying space of a reductive group was introduced by Totaro in [Tot], where he calculated the Chow rings of the classifying spaces of several finite groups and algebraic groups including $O(n)$, $Sp(2n)$, etc. Edidin and Graham in [E-G] introduced the equivariant Chow ring. Rojas and Vistoli in [R-V], using the techniques of equivariant Chow groups, calculated the Chow ring $A^*_{SO(n)}$ in case $n$ is even (the odd case was already addressed in Pandharipande [Pan] and Totaro [Tot]). 

Holla and Nitsure in [H-N] considered the general orthogonal group $GO(n, \mathbb C)$, which is also called the group of similitudes, and calculated the singular cohomology ring of its classifying space $H^*_{\rm sing}(BGO(n,\mathbb C);\ftwo)$. In [Bh], the author of the present paper considered the algebraic version of the above Lie group, namely the general orthogonal group scheme $GO(n)$ over an algebraically closed field of characteristic different from $2$, and showed that the smooth-\'etale cohomology ring $H^*_{\rm sm\mhyphen et}(BGO(n);\ftwo)$ of the algebraic stack $BGO(n)$ has the same description in terms of generators and relations over $\ftwo$ as the singular cohomology ring computed by Holla and Nitsure in [H-N]. In this present note, we calculate the Chow ring of the classifying space of $GO(n)$ over a field of characteristic different from $2$ in the sense of Totaro [Tot], using the methods of equivariant Chow groups. By the results of Totaro [Tot], this ring can be canonically identified with the ring of characteristic classes for principal $GO(n)$-bundles on smooth, quasi-projective schemes. In other words, the Chow ring of the classifying space of $GO(n)$ is the ring of all intersection theoretical invariants for families of line bundle valued nondegenerate quadratic forms.

Henceforth, all schemes and morphisms are over a fixed field $k$ (not necessarily algebraically closed) with characteristic different from $2$. We recall the definition of the algebraic group $GO(n)$ over $k$. Let $V=k^n$, and let $q:V\to k$ be the quadratic form, defined by \[q(x_1,\ldots,x_{2m})=x_1x_{m+1}+\cdots+x_mx_{2m},\]for the even case $n=2m$, and by \[q(x_1,\ldots,x_{2m+1})=x_1x_{m+1}+\cdots+x_mx_{2m}+x_{2m+1}^2,\] for the odd case $n=2m+1$. Let $GO(n)$ be the affine algebraic group scheme of invertible linear automorphisms of $V$ that preserve the quadratic form $q$ up to a scalar. In terms of matrices, let $J$ denote the nonsingular symmetric matrix of the bilinear form corresponding to $q$. Then as a functor of points, $GO(n)$ attaches to each $k$-algebra $S$ the group 
\[GO(n)(S)=\{A\in GL_n(S)\colon \exists\, a\in S^\times, \; {}^t\!AJA=aJ\}.\]  The algebraic group $GO(n)$ is reductive, since its defining representation on $k^n$ is irreducible. Note that if $k'/k$ is a field extension such that the quadratic form $q$ extended to $V\otimes_k k'=k'^{n}$ is equivalent to the quadratic form $\sum_i x_i^2$, given by the identity matrix $I_n$, then over $k'$, the algebraic group $GO(n)$ defined above is isomorphic to the algebraic group $GO(n)$ defined in [Bh].

The scalar $a$ in the definition determines the character $\sigma:GO(n)\to \gm$ that satisfies ${}^t\!AJA=\sigma(A)J$. Given a scheme $X$, and a principal $GO(n)$-bundle $P$ on $X$ (locally trivial in the \'etale topology), consider the rank $n$ vector bundle $E$ associated to the defining representation $GO(n)\subset GL_n$, and the line bundle $L$ determined by the character $\sigma$. The nondegenerate symmetric bilinear form corresponding to $q$ induces a nondegenerate symmetric bilinear form $b:E\otimes_{\mathcal O_X} E\to L$. Conversely, given a \emph{nondegenerate quadratic triple of rank $n$} $(E,L,b)$, which is a triple consisting of a vector bundle $E$ of rank $n$, a line bundle $L$ and a nondegenerate symmetric bilinear $b:E\otimes_{\mathcal O_X} E\to L$, we can reduce the structure group of $E$ to $GO_n$ by applying Gram-Schmidt orthonormalization \'etale locally on $X$. 

Let $A^*_G$ denote the Chow ring of the classifying space of a reductive group $G$ in the sense of Totaro [Tot]. Note that for any $n\ge 1$, there is a canonical isomorphism 
\[S\!O(2n+1)\times \gm\isomor GO(2n+1).\] Note that $B\gm$ is approximated by the projective spaces $\mathbb P^m_k$ in the sense of Totaro [Tot], and we have for any smooth scheme $X$ the following natural isomorphisms.
\[A^*(X)\otimes A^*(\mathbb P^m_k)\isomor A^*(X\times \mathbb P^m_k)\]Since $BSO(2n+1)$ is approximated by smooth schemes, there is the K\"unneth isomorphism 
\[A^*_{SO(2n+1)}\otimes A^*_{\gm}\isomor A^*_{GO(2n+1)}.\]
This determines the Chow ring for $GO(2n+1)$, because $A^*_{\gm}\cong \mathbb Z[\lambda]$, and the Chow ring for $SO(2n+1)$ is given by [Pan] and [Tot]. Therefore we are left with the task of calculating the Chow ring only in the even case $GO(2n)$. The rest of this note is devoted to the calculation of $A^*_{GO(2n)}$.

In Section 2, we recall Totaro's definition of the Chow ring of a classifying space from [Tot], and the basic notions of equivariant Chow groups from Rojas and Vistoli [R-V]. In Section 3, we calculate $A^*_{GO(2n)}$ in terms of explicit generators and relations over $\mathbb Z$. We show, in terms of quadratic triples $(E,L,b)$, $A^*_{GO(2n)}$ is generated by the Chern classes $c_i(E)$ and the Chern class $\lambda$ of $L$. The nondegenerate symmetric bilinear form determines an isomorphism of vector bundles $E\stackrel\sim\to E^\vee\otimes L$. This isomorphism gives the following relations among Chern classes of $E$ and the class $\lambda$. \[c_p=\sum_{i=0}^p(-1)^i{2n-i\choose p-i}c_i \lambda^{p-i},\;p=1,\ldots,2n\] The main theorem of this note is Theorem \ref{maintheorem}, which asserts that the ideal of relations in $A^*_{GO(2n)}$ is generated by the above relations. 

When the base field is $\mathbb C$, the invariants in $A^*_{GO(2n)}$ transform to the corresponding classes in $H^*\hsing(BGO(2n),\ftwo)$ under the cycle map (see Remark \ref{cycleclass}).

To prove the main theorem (Theorem \ref{maintheorem}), we first observe that $\lambda$ and the even Chern classes are algebraically independent in $A^*_{GO(2n)}$. Eventually we focus on the torsion subgroup of $A^*_{GO(2n)}$ (which is an ideal) in order to carry out the rest of the proof. The inclusion $O(2n)\subset GO(2n)$, where $O(2n)$ is looked upon as the algebraic group of linear automorphisms of $k^{2n}$ that preserve the quadratic form $\sum_i x_ix_{i+n}$, gives rise to a ring homomorphism $A^*_{GO(2n)}\to A^*_{O(2n)}$. A biproduct of the proof of the main theorem is that this homomorphism determines an isomorphism of the corresponding torsion subgroups.
\vskip 1em
\noindent{\bf Acknowledgments.} The author would like to thank Yogish Holla and Nitin Nitsure for useful discussions and suggestions. 

\section{Basic notions recalled}
\subsection{Chow ring of the classifying space}
Totaro [Tot] defined the Chow groups of the classifying space of a reductive algebraic group $G$ as follows. Let $N>0$ be an integer. Then there is a finite dimensional linear representation $V$ of $G$, with a $G$-equivariant closed subset $S$ of codimension $\ge N$ such that the action of $G$ on $(V-S)$ is free, such that the quotient $(V-S)/G$ exists in the category of schemes and $(V-S)\to (V-S)/G$ is a principal bundle, locally trivial in \'etale topology. For any such other pair $(V',S')$, there are canonical isomorphisms for $i<N$
\[A^i\left (\frac{(V-S)}{G}\right)\isomor A^i\left (\frac{(V'-S')}{G}\right).\]The $i$-th Chow group $A^i_G$ of the classifying space of $G$ is defined to be $A^i\left (\frac{V-S}{G}\right)$, where $(V,S)$ satisfies the above condition. The isomorphism above shows that $A^i_G$ is independent of the choice of the particular pair $(V,S)$. The graded group $A^*_G$ naturally has the structure of a graded ring under intersection product.

Note that in the above approach, one defines the Chow ring of the classifying space without referring to a possible classifying space $BG$ whether in the category of algebraic stacks or simplicial spaces.

\subsection{Equivariant Chow groups}
Let $X$ be a scheme on which a reductive algebraic group $G$ acts. Suppose $V$ is a finite dimensional linear representation of $G$, and $S\subset V$ is a $G$-equivariant closed subset such that the induced action of $G$ on $(V-S)$ is free, with a principal bundle quotient $(V-S)\to (V-S)/G$. If the codimension of $S$ in $V$ is $N$, then for $i<N$, the equivariant Chow groups of $X$ are defined as 
\[A^i_G(X):=A^i\left (\frac{X\times (V-S)}{G}\right).\]This definition is independent of the particular choices of $V$ and $S$.

The Chow ring of the classifying space is recovered by taking $X$ to be the point $\spec k$, as $A^i_G(\spec k)=A^i_G$.

Suppose $G$ fits into an exact sequence of reductive algebraic groups as follows
\[1\to H\to G\stackrel\chi\to \gm\to 1.\] Consider the action of $G$ on $\mathbb A^1$ by the character $\chi$. The localisation sequence for $\gm\subset \mathbb A^1$ is \[A^*_G(\spec k)\to A^*_G(\mathbb A^1)\to A^*_G(\gm)\to 0.\] The quotient $(\gm\times (\mathbb A^1-\{0\}))/G$ can be naturally identified with $(\mathbb A^1-\{0\})/H$. Therefore the above localization sequence gives an exact sequence
\[A^*_G\stackrel{c}\to A^*_G\to A^*_H\to 0,\]where $c$ is the Chern class of the line bundle given by the character $\chi$.

\section{Calculation of $A^*_{GO(2n)}$}
\centerline{\large\bf Generators for $A^*_{GO(2n)}$}
\vskip 1em
Recall that we have a short exact sequence of reductive algebraic groups 
\[1\to O(2n)\to GO(2n)\stackrel\sigma\to \gm\to 1,\]which, by the results of the last section, gives an exact sequence
\begin{eqnarray}\label{localization}A^*_{GO(2n)}\stackrel{\lambda}\to A^*_{GO(2n)}\to A^*_{O(2n)}\to 0\end{eqnarray}where $\lambda$ is the Chern class corresponding to the character $\sigma$.

Rojas and Vistoli [R-V] showed that \[A^*_{O(2n)}\cong \frac{\mathbb Z[c_1,\ldots,c_{2n}]}{(2c_{\rm odd})},\] hence $A^*_{O(2n)}$ is generated by the Chern classes $c_1,\ldots,c_{2n}$. The exact sequence (\ref{localization}) shows that $A^*_{GO(2n)}$ is generated by $c_1,\ldots,c_{2n},\lambda$.

\vskip 1em \centerline{\large\bf Relations among $c_1,\ldots,c_{2n},\lambda$}
\vskip 1em

Recall that to give a principal $GO(2n)$-bundle is to give a triple $(E,L,b)$ consisting of a vector bundle $E$ of rank $2n$, a line bundle $L$ and a nondegenerate symmetric form $b:E\otimes E\to L$. The form $b:E\otimes E\to L$ determines an isomorphism of vector bundles $E\isomor E^\vee\otimes L$, which gives the following relations among $\lambda$ and the Chern classes of $E$.
\begin{eqnarray}\label{relation}c_p&=&\sum_{i=0}^p(-1)^i{2n-i\choose p-i}c_i \lambda^{p-i},\;p=1,\ldots,2n\end{eqnarray}

Let $R$ be the quotient of the polynomial ring $\mathbb Z[\lambda,c_1,\ldots, c_{2n}]$ by the ideal generated by the above $2n$ relations. Let $q:R\to A^*_{GO(2n)}$ be the ring homomorphism that sends $\lambda$ to the Chern class of the line bundle defined by the character $\sigma$, and $c_i$ to the $i$th Chern class of the defining representation $GO(2n)\subset GL_{2n}$.

\begin{theorem}\label{maintheorem} The map $q:R\to A^*_{GO(2n)}$ is an isomorphism.\end{theorem}
Before giving the proof of this theorem, we make a remark.
\begin{remark}\label{cycleclass}
For a nonsingular variety $X$ over the field of complex numbers, the cycle map is a homomorphism of graded rings $cl^X:A^*(X)\to H^{2*}\!\hsing(X,\mathbb Z)$, functorial in $X$ (see Fulton [Fu], Chapter 19). By composing with the change of coefficients map $H^{2*}\!\hsing(X,\mathbb Z)\to H^{2*}\!\hsing(X,\ftwo)$, we get a homomorphism of graded rings $A^*(X)\to H^{2*}\!\hsing(X,\ftwo)$, which we will denote by $\bar {cl}^X$. If $E$ is a vector bundle on $X$, then $\bar {cl}^X(c_i(E))=\bar c_i(E)\in H^{2i}\!\hsing(X,\ftwo)$, the mod-2 reduced Chern classes. For a reductive group $G$, this gives rise to a homomorphism of graded rings (see [Tot]), which is functorial in the group $G$
\[\bar {cl}^G:A^*_G\to H^{2*}\!\hsing(BG,\ftwo).\] For $G=GO(2n,\CC)$, the Chern classes $c_i\in A^*_{GO(2n)}$ transform under the above cycle map to the corresponding classes in $H^{2*}(BGO(2n,\mathbb C),\ftwo)$ described by Holla and Nitsure (see Section 3 of [H-N-2]).

\end{remark}

{\footnotesize{\bf Note }In case of \'etale cohomology, I am ignorant whether the cycle map $Z^i(X)\to H^{2i}\het(X,\ftwo)$ passes through rational equivalence, although it is due to Grothendieck that the cycle map $Z^i(X)\to H^{2i}(X,\mathbb Q_\ell)$ does pass through rational equivalence. Therefore I do not have the means to compare the Chern classes $c_i\in A^i_{GO(2n)}$ with the the images of $\bar c_i$ in $H^{2i}\hliset(BGO(2n),\ftwo)$ which have the similar description as in Proposition 3.2 of [H-N-2] (see [Bh] Section 5).
}

\vskip 1em
The rest of this article is devoted to the proof of the Theorem \ref{maintheorem}, which is the main result of this article. Let us begin by recalling that the map $q$ is surjective. We will prove that it is injective. The plan of the proof is as follows. We will first show that the elements $\lambda, c_2,c_4,\ldots, c_{2n}$ are algebraically independent in $A^*_{GO(2n)}$ (Corollary \ref{indep}). Then we prove Lemma \ref{bdens} and Corollary \ref{maincorol}, which will imply that it is enough to prove the injectivity only for the torsion part. We will complete the proof of injectivity for the torsion part in a few steps. As a concluding remark, we observe (Corollary \ref{ftwo}) that the torsion subgroup is in fact an $\ftwo$-vector space.

\begin{remark}\label{twomult} Note that for each \emph{odd} $p$, we have the following identity in $R$ \begin{eqnarray}\label{odd}2c_p=\sum_{i=0}^{p-1}(-1)^ic_i{2n-i\choose p-i} \lambda^{p-i}\end{eqnarray}In particular, $(2c_{odd})R\subset \lambda R$.
\end{remark}
\begin{lemma}Let the free polynomial algebra $B=\mathbb Z[\lambda,c_2,c_4,\ldots,c_{2n}]$ be given the grading where $\lambda$ has homogeneous degree $1$ and each $c_{2i}$ has homogeneous degree $2i$. Let $D$ be a graded domain, and let $\phi:B\to D$ be a graded homomorphism such that \\
(1) the restriction of $\phi$ to $\mathbb Z[c_2,c_4,\ldots,c_{2n}]$ is an injection, \\
(2) so is the composite $\mathbb Z[c_2,c_4,\ldots,c_{2n}]\to D\to D/\phi(\lambda)D$, and \\
(3) $\phi(\lambda)\in D$ is non-zero.\\
Then $\phi$ is an injection.\end{lemma}
\proof For a polynomial in $B$ with degree $\le 1$, the image is always non-zero in $D$, by (1) and (3) of the hypothesis. We will prove the injectivity of $\phi$ by induction on the degree, as we know that the injectivity holds in degree $\le 1$.

Suppose the injectivity holds for degree $<r$, and suppose $f\in B$ is a homogeneous polynomial of degree $r$ with $\phi(f)=0$. We can write $f=\lambda\cdot g+h$, where $g\in B$ in homogeneous of degree $(r-1)$, and $h\in\mathbb Z[c_2,c_4,\ldots,c_{2n}]$ is of degree $r$. By (2), $h=0$. Therefore, $f=\lambda\cdot g$. Suppose $g\ne 0$. Since ${\rm deg}(g)<r$, by induction hypothesis, $\phi(g)\ne 0$. But since $D$ is a domain, and as by (3) $\phi(\lambda)\ne 0$, we see that $\phi(f)=\phi(\lambda)\phi(g)\ne 0$, a contradiction, so that $g=0$.\hfill$\Box$

Let us consider the invertible $2n\times 2n$ matrix $J=(0,I_n\,;I_n,0)$. By definition, $GO(2n)$ is the group scheme of invertible matrices $A$ such that ${}^t\!AJA=aJ$ for some scalar $a$. Consider the closed subgroup scheme $\Gamma\subset GO(2n)$, consisting of matrices of the following form \[\left(\begin{array}{llllll}\lambda t_1 &&&&&\\ &\ddots&&&&\\ &&\lambda t_n&&& \\ &&&t_1^{-1}&&\\ &&&&\ddots &\\ &&&&& t_n^{-1}\end{array}\right)\]Then $\Gamma$ is a maximal torus in $GO(2n)$. The inclusion $\Gamma\subset GO(2n)$ induces a homomorphism of the corresponding Chow rings \[A^*_{GO(2n)}\to A^*_\Gamma=\mathbb Z[\lambda,t_1,\ldots,t_n],\] sending $\lambda$ to $\lambda$. Applying the above lemma to $D=A^*_\Gamma$ and the graded homomorphism which is the composite $B\to A^*_{GO(2n)}\to D$, and using the fact that $c_2,\ldots,c_{2n}$ are algebraically independent in the ring $A^*_{O(2n)}\cong A^*_{GO(2n)}/(\lambda)$, we get the following.

\begin{corol}\label{indep}The elements $\lambda,c_2,c_4,\ldots,c_{2n}$ are algebraically independent in the ring $A^*_{GO(2n)}$.\end{corol}

\noindent{\bf Notation.} Let us recall that we adopted the notation $R$ for $\mathbb Z[\lambda,c_1,\ldots,c_{2n}]/I$, where $I$ is generated by the relations coming from $E\isomor E^\vee\otimes L$ enlisted in equation (\ref{relation}). We have the surjective map $q:R\to A^*_{GO(2n)}$. Let $A$ denote $A^*_{GO(2n)}$ and let $B$ denote the polynomial ring $\mathbb Z[\lambda,c_{even}]$. Since the composite $B\to R\to A$ is injective by the corollary, the map $B\to R$ is also injective. Hence we might consider $B\subset R$ as well as $B\subset A$. We will denote the torsion ideals of $R$ and $A$ by $\tr$ and $\ta$, respectively.

\begin{lemma}\label{bdens}For any $a\in R$, there is some $s>0$ such that $2^sa\in B$. Similarly, for any $a\in A$, there is some $s>0$ such that $2^sa\in B$.\end{lemma}
\proof We show this by induction on the degree of $a$. It is enough to prove the lemma for all odd Chern classes i.e. for $a=c_{2i+1}$. By equation (\ref{odd}), $2c_{2i+1}=g_i(\lambda,c_1,\ldots,c_{2i})$, a polynomial in $\lambda$ and lower degree Chern classes. Now by the induction hypothesis, for each $j<i$, there is some $m_j> 0$ such that $2^{m_j}c_{2j+1}\in A$. If $n_j$ is the greatest index with which $c_{2j+1}$ occurs in $g_i$, then $s=\sum_jn_jm_j$ serves the purpose.\hfill$\Box$

We make the following observations as corollaries.
\begin{corol}\label{maincorol}
(a) The groups $\tr$ and $\ta$ consist of $2$-primary elements only. 

(b)  We have $\ker(q)\subset \tr$.

(c) We have $q\inv(T_A)=T_R$. In particular, we have the induced isomorphism of free abelian group \[R/T_R\cong A/T_A.\] 

(d) If $a\in R$ such that $\lambda a\in \tr$, then $a\in\tr$. Similarly, if $a\in R$ such that $\lambda a\in \tr$, then $a\in\tr$.

(e) $\ker(q)\subset \tr\cap \lambda R$.

(f) We have $\tr=\ker(R\to A^*_\Gamma))$ and $\ta=\ker(A\to A^*_\Gamma)$.

(g) There are isomorphisms $R/(\lambda)\isomor A/(\lambda)\isomor A^*_{O(2n)}$.
\end{corol}

\proof

(a) Follows from the lemma.

(b) Suppose $a\in R$ is such that $q(a)=0$. Then there is some $m\ge 0$ such that $0=2^mq(a)=q(2^ma)$, and $2^ma\in B$. Hence $2^ma=0$. This shows that the kernel of $q$ has only $2$-primary elements. 

(c) Follows from (b). For, if $a\in q\inv(\ta)$, then $2^sa\in\ker(q)\subset\tr$, so $a\in\tr$.

(d) Otherwise, by the lemma, we get $2^sa\in B$ and $2^sa\ne 0$. This means $\lambda a$ is never a torsion, because $2^s\lambda a\in B-\{0\}$.

(e) To see this, first note that any element $a$ of $\mathbb Z[\lambda,c_1,\ldots,c_{2n}]$ is written as $\alpha+\lambda \beta+\gamma$, where $a\in\mathbb Z[c_{even}]$, and $\gamma\in (c_{odd})\mathbb Z[c_1,\ldots,c_{2n}]$. If $\alpha\ne 0$, then see that the image of $2a$ in $A^*_{O(2n)}$ is equal to the image of $2\alpha$, which is non-zero. Therefore if $a\in\ker(q)$, then $\alpha=0$ and $a=\lambda\beta+\gamma$. Now, the image of $a$ in $A^*_{O(2n)}$ vanishes, so $\gamma\in (2c_{odd})\mathbb Z[c_1,\ldots,c_{2n}]$. But by equation (\ref{odd}), we have $(2c_{odd})R\subset \lambda R$. So the assertion follows.

(f) That $\tr\subset \ker(R\to A^*_\Gamma)$ is obvious, because $A^*_\Gamma$ has no torsions. To see the other inclusion, let $a\in \ker(R\to A^*_\Gamma)$. If $a$ was not a torsion, then by lemma, there is $s>0$ such that $2^sa\in B-\{0\}$. But $B\to A^*_\Gamma$ is injective as we have already seen. Similarly for $\ta$.

(g) The map $R/(\lambda)\to A^*_{O(2n)}$ is an isomorphism by the definition of $R$. On the other hand, equation (\ref{localization}) shows that $A/(\lambda)\to A^*_{O(2n)}$ is an isomorphism. \hfill $\Box$
\vskip 1em

{\bf Proof of the main theorem.} 

We only have to prove that $\tr\cap\lambda R=0$, since $\ker(q)\subset \tr\cap \lambda R$.

{\bf Step 1.} Let $\lambda_A$ denote the multiplication by $\lambda:A\to A$. In what follows, $\ker\lambda$ will denote the kernel of the multiplication $\lambda:R\to R$, while $\ker\lambda_A$ will denote the kernel of $\lambda_A:A\to A$. We have the following short exact sequences, where the right side maps are multiplications by $\lambda$
\[0\to \ker\lambda\cap\tr\to \tr\to\lambda R\cap \tr\to 0,\]
\[0\to \ker\lambda_A\cap\ta\to \ta\to\lambda A\cap \ta\to 0.\]

Indeed, we need only see that the right side map is surjective. But if $x\lambda\in\tr$, then $x$ has to be a torsion element by (d) of Corollary \ref{maincorol}. Similar reasons apply to the latter sequence.

{\bf Step 2.} If $C$ denotes the image in $A^*_\Gamma$ of the composite $R\stackrel q\to A\to A^*_\Gamma$ (which we will call $\pi$), then we have the following short exact sequence by (f) of Corollary \ref{maincorol}
\[0\to \tr\to R\stackrel\pi\to C\to 0.\]
Now, $C\subset A^*_\Gamma$, so that multiplication by $\lambda$ is injective on $C$. Therefore, if $x\in R$ such that $\lambda x=0$, then $\pi(x)$ cannot be nonzero in $C$. Therefore, we get the following two inclusions, of which the latter follows by a similar argument \[\ker\lambda\subset \tr\]
\[\ker\lambda_A\subset\ta.\]

{\bf Step 3.} Let $T_O$ denote the subgroup of all torsion elements in $A^*_{O(2n)}$. Under the composite $R\to R/\lambda R\isomor A^*_{O(2n)}$, torsion elements map inside $T_O$. This gives a map $\tr\to T_O$. We will shortly show that this is surjective. As a consequence, we will have a short exact sequence \[0\to \lambda R\cap \tr\to \tr\to T_O\to 0.\]Again, the surjectivity of the composite $\tr\to\ta\to T_O$ will imply the surjectivity of $\ta\to T_O$, which, together with the isomorphism $A/\lambda\isomor A^*_{O(2n)}$ will give another short exact sequence
 \[0\to\lambda A\cap\ta\to \ta\to T_O\to 0.\]

Now let us go back to the proof of the surjectivity of $\tr\to T_O$. Since we have $T_O\isomor (c_{odd})A^*_{O(2n)}$ as graded groups, it is sufficient to show that for each odd $p<2n$, there is a torsion element $\beta_p$ in $R$ such that $\beta_p\mapsto c_p$ under $\tr\to T_O$. From equation (\ref{relation}) that in $R$, for each odd $p$, we have the following equality
\begin{eqnarray}\label{step3}c_{p+1}=c_{p+1}-(2n-p)c_p\lambda +\lambda^2\alpha'_p\end{eqnarray}for some $\alpha'_p\in R$, so that $\lambda((2n-p)c_p-\lambda\alpha'_p)=0$. But since $\ker\lambda\subset\tr$, we see that the element $\beta_p=(2n-p)c_p-\lambda\alpha'_p$ is torsion. Since $(2n-p)$ is odd, and since $2c_p=0$ in $T_O$, we see that $\beta_p\mapsto c_p$ under $\tr\to T_O$, as desired.

{\bf Step 4.} We have these two short exact sequences, which come from Step 1, by the substitutions $\ker\lambda=\ker\lambda\cap\tr$ and $\ker\lambda_A=\ker\lambda_A\cap \ta$.

\[0\to \ker\lambda\to \tr\to \lambda R\cap\tr\to 0\]
\[0\to \ker\lambda_A\to \ta\to \lambda A\cap \ta\to 0\]

{\bf Step 5.} Note that both $R$ and $A$ are noetherian rings, and their ideals $\tr$ and $\ta$ are finitely generated graded ideals. By (a) of \ref{maincorol} torsion elements are $2$-primary. So there is some $N>0$ such that $2^N\tr=0$, and $2^N\ta=0$. For each $m$, the graded pieces $R_m$ and therefore $A_m$ are finitely generated abelian groups. Therefore for each $m$, the graded pieces $(\tr)_m$ and $(\ta)_m$, which are finitely generated abelian group and hence finitely generated $\mathbb Z/2^N\mathbb Z$-modules, are finite sets. As $\ker\lambda\subset\tr$ and $\ker\lambda_A\subset\ta$, the sets $(\ker\lambda)_m$ and $(\ker\lambda_A)_m$ are finite as well. The group $T_O$ is actually an $\ftwo$-vector space and, for similar reasons as above, each $(T_O)_m$ is a finite set.

Therefore, from the two short exact sequences listed in Step 3, we get
\[\#(\tr)_m-\#(\lambda R\cap \tr)_m=\#(T_O)_m=\#(\ta)_m-\#(\lambda A\cap \ta).\]

From those listed in Step 4, we get
\[\#(\tr)_m-\#(\lambda R\cap \tr)_m=\#(\ker\lambda)_m,\]
\[\#(\ta)_m-\#(\lambda A\cap \ta)_m=\#(\ker\lambda_A)_m.\]

Therefore for each $m$, \[\#(\ker\lambda)_m=\#(\ker\lambda_A)_m=\#(T_O)_m.\] 

{\bf Step 6.} We finally proceed to prove that $R\to A$ is injective (therefore bijective). By induction on degree, we will prove that $R_m\to A_m$ is injective (therefore bijective) for each $m$. 

To begin the induction, note that this is true for $m=0$ and $m=1$, by \ref{indep}. Now, suppose this is true for $1,\ldots,m$, and we will prove it for $m+1$.

Since $R_m\to A_m$ is bijective by assumption, $(\ker\lambda)_m\to (\ker\lambda_A)_m$ is injective. But they have the same number of elements by Step 5. Hence  $(\ker\lambda)_m\to (\ker\lambda_A)_m$ is surjective (therefore bijective).

With this, and the fact that $R/\lambda R\isomor A/\lambda A$, we have the following commutative diagram, whose left and right side vertical maps are isomorphisms, and whose rows are exact.

\centerline{\xymatrix{
0\ar[r]&(\ker\lambda)_m\ar[r]\ar[d]_\cong& R_m\ar[d]_{q}^\cong\ar[r]^{\lambda}& R_{m+1}\ar[d]_{q}\ar[r] & (R/\lambda R)_m\ar[d]^{\cong}\ar[r]&0\\
0\ar[r]&(\ker\lambda_A)_m\ar[r] & A_m\ar[r]_\lambda& A_{m+1}\ar[r]& (A/\lambda A)_m\ar[r]&0
}}

By five lemma, $q:R_{m+1}\to A_{m+1}$ is an isomorphism. \hfill$\Box$

\begin{remark}\label{lambda} For each odd number $2p+1$, there is an odd number $n_{2p+1}$ such that $n_{2p+1}\lambda c_{2p+1}=\lambda f_p(c_{2p},c_{2p-1},\ldots,c_1,\lambda)$, where $f_p$ is a polynomial. Indeed, in equation (\ref{step3}) in Step 3 of the proof of the main theorem, the term $\alpha'_p$ is a polynomial in lower dimensional Chern classes and $\lambda$, so one can do induction on $p$. Now, similarly, $n_{2p-1}n_{2p+1}\lambda c_p=\lambda f'_p$, where $f'_p$ is a polynomial of $c_{2p},c_{2p-2},c_{2p-3},\ldots,c_1,\lambda$. In this way, there is an odd number $N_{2p+1}$ such that $N_{2p+1}\lambda c_{2p+1}=\lambda g_{2p+1}$, where $g_{2p+1}\in B$. Therefore we can say that if $\gamma=\lambda\gamma'\in A$, then there is some odd $N''$ such that $N''\gamma\in B$. \end{remark}

\begin{lemma} Given any $\gamma\in A$, there is an odd number $N$ such that $2N\gamma\in B$.\end{lemma}
\proof % We will prove this lemma by induction on the highest odd $m$ such that $c_m$ appears in the expression of $\gamma$.
To prove our lemma, is enough to assume that $\gamma$ can be given by a monomial, involving some odd Chern classes. So $\gamma=c_{2i_1+1}^{m_1}\ldots c_{2i_s+1}^{m_s}g$, where $g$ does not involve any $c_{\rm odd}$. It is also enough to assume that $g\equiv 1$. Now, 
\[2\gamma= f_{i_1}(c_{2i_1},\ldots, c_1,\lambda) \lambda c_{2i_1+1}^{m_1-1}\ldots c_{2i_s+1}^{m_s},\]where $f$ is a polynomial. Hence there is some odd number $N'$ such that \[2N'\gamma=\lambda f_{i_1}(c_{2i_1},\ldots, c_1,\lambda)g',\] where $g'\in B$. Now, by the remark preceding our lemma, there is some odd number $N''$ such that $2N'N''\gamma\in B$. \hfill$\Box$
\begin{corol}\label{ftwo} The torsion subgroup $\ta\subset A^*_{GO(2n)}$ is an $\ftwo$-vector space.\end{corol}
\proof This follows from the last lemma and the fact that each element in $\ta$ is $2$-primary.\hfill$\Box$

\begin{remark} Remark \ref{lambda} also shows that $\lambda T_A=\lambda R\cap\ta=0$. Indeed, from \ref{lambda}, if $\gamma\in\lambda \ta$, then there is an odd number $N$ such that $N\gamma=0$. But $\gamma\in\ta$ also, hence is a $2$-torsion. So $\gamma=0$, as desired. Consequently, from Step 3 and Step 4 of the proof of the main theorem,\[\ker\lambda=\ta\isomor T_O.\]\end{remark}

\section*{References}
\noindent [Bh] Saurav Bhaumik : \emph{Characteristic classes for GO(2n) in \'etale cohomology}, Proc. Indian Acad. Sci. (Math. Sci.), Volume 123, No. 2 (May 2013), pp 225-233

\noindent [E-G] Dan Edidin and William Graham : \emph{Equivariant intersection theory}, Invent. Math. 131 (1998), no. 3, 595-634

\noindent [Fu] William Fulton : \emph{Intersection Theory}, Second Edition 1998, Springer

\noindent [H-N] Yogish I Holla and Nitin Nitsure : \emph{Characteristic Classes for $GO(2n,\mathbb C)$}, Asian J. Math. 5 (2001) 169-182

\noindent [H-N-2] Holla, Y.I. and Nitsure, N. : Topology of quadric bundles. Internat. J.
Math. 12 (2001), 1005-1047.

\noindent [Pan] Rahul Pandharipande : \emph{Equivariant Chow rings of O(k), SO(2k+1), and SO(4)},
J. Reine Angew. Math. 496 (1998), 131-148.

\noindent [Tot] Burt Totaro : \emph{Chow Ring of a Classifying Space}, Algebraic $K$-Theory, Proc. Sympos. Pure Math., 67, Amer. Math. Soc. (1999).

\noindent [R-V] Luis Alberto Molina Rojas and Angelo Vistoli : \emph{On the Chow Rings of Classifying Spaces for Classical Groups}, Rend. Sem. Mat. Univ. Padova, Vol. 116 (2006)
\end{document}